# Weighted Sum Rate Optimization for Downlink Multiuser MIMO Systems with per Antenna Power Constraint: Downlink-Uplink Duality Approach

Tadilo Endeshaw Bogale and Luc Vandendorpe, Email: {tadilo.bogale, luc.vandendorpe}@uclouvain.be
ICTEAM Institute, Universitè catholique de Louvain, Place du Levant B-1348, Louvain La Neuve, Belgium

*Abstract*— This paper considers weighted sum rate maximization constrained with a per base station (BS) antenna power problem for multiuser multiple-input multiple-output (MIMO) systems. For this problem, we propose new downlink-uplink duality based solution. We solve the problem as follows. First, by introducing additional optimization variables, we reformulate our problem into an equivalent problem that incorporates a weighted sum mean-square-error (MSE) term. Second, we establish novel weighted sum MSE duality. The duality is established by modifying the input covariance matrix of the dual uplink problem, and formulating the noise covariance matrix of the uplink channel as a fixed point function. Third, we optimize the introduced variables and powers in the downlink channel by a Geometric Program (GP) method. Fourth, using the duality result and the solution of GP, we apply alternating optimization technique to solve the original downlink problem. In our simulation results, we have observed that the proposed duality based solution utilizes less power than that of existing algorithms.

## I. INTRODUCTION

Multiple-input multiple-output (MIMO) system is a promising approach to fully exploit the spectral efficiency of wireless channels. The capacity region of a MIMO broadcast channel is achieved by dirty paper coding (DPC) technique of [1]. However, since DPC is non-linear, practical realization of DPC is difficult. Due to this, linear processing is motivated as it exhibits good performance/complexity trade-off.

In [2] and [3], weighted sum rate maximization constrained with a per-BS antenna power problem has been examined. Furthermore, simulation results demonstrate that the approach of the latter paper has better achievable weighted sum rate than that of [2]. Both of these papers examine their problems directly in the downlink channel. In the downlink channel, since the precoders of all users are jointly coupled, downlink transceiver design problems have more complicated mathematical structure than that of the dual uplink problems [4], [5]. In [6], we have demonstrated that for a per BS antenna power constrained sum MSE-based problems, the duality approach of [6] requires less total BS power than the approach of [7]. Due to these reasons, duality approach for solving downlink problems has received a lot of attention.

In this paper, we examine the weighted sum rate maximization constrained with a per-BS antenna power problem (P1). For this problem, we propose new duality based

The authors would like to thank the Region Wallonne for the financial support of the project MIMOCOM in the framework of which this work has been achieved.

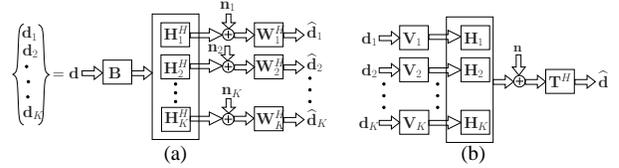

Fig. 1. Multiuser MIMO system. (a) downlink channel. (b) uplink channel.

iterative solution. Our duality based iterative solution for $\mathcal{P}1$ is explained as follows. First, by introducing additional optimization variables, we reformulate $\mathcal{P}1$ into an equivalent problem that incorporates a weighted sum MSE term like in [3]. Second, we establish novel weighted sum MSE duality. The duality is established by modifying the input covariance matrix of the dual uplink problem, and formulating the noise covariance matrix of the uplink channel as a fixed point function. Third, we optimize the introduced variables and powers in the downlink channel by a GP method. Fourth, using the duality result and solution of GP, we utilize alternating optimization technique to solve the original downlink problem. In our simulation results, we have demonstrated that the proposed duality based algorithm utilizes less total BS power than that of existing algorithms.

## II. SYSTEM MODEL

In this section, multiuser MIMO downlink and virtual uplink system models are considered. The BS and $k$th mobile station (MS) are equipped with $N$ and $M_k$ antennas, respectively. The total number of MS antennas are thus $M = \sum_{k=1}^{K} M_k$. By denoting the symbol intended for the $k$th user as $\mathbf{d}_k \in \mathbb{C}^{S_k \times 1}$ and $S = \sum_{k=1}^{K} S_k$, the entire symbol can be written in a data vector $\mathbf{d} \in \mathbb{C}^{S \times 1}$ as $\mathbf{d} = [\mathbf{d}_1^T, \cdots, \mathbf{d}_K^T]^T$. The BS precodes $\mathbf{d}$ into an $N$ length vector by using its overall precoder matrix $\mathbf{B} = [\mathbf{b}_{11}, \cdots, \mathbf{b}_{KS_K}]$, where $\mathbf{b}_{ki} \in \mathbb{C}^{N \times 1}$ is the precoder vector of the BS for the $k$th MS $i$th symbol. The $k$th MS employs a receiver $\mathbf{w}_{ki}$ to estimate its symbol $d_{ki}$. We follow the same channel matrix notations as in [6]. The estimate of the $k$th MS $i$th symbol ($\hat{d}_{ki}$) is given by

$$\hat{d}_{ki} = \mathbf{w}_{ki}^H (\mathbf{H}_k^H \mathbf{B} \mathbf{d} + \mathbf{n}_k) \qquad (1)$$

where $\mathbf{H}_k^H \in \mathbb{C}^{M_k \times N}$ is the channel matrix between the BS and $k$th MS, and $\mathbf{n}_k$ is the additive noise at the $k$th MS. We assume that the entries of $\mathbf{d}_k$ are i.i.d zero-mean circularly symmetric complex Gaussian (ZMCSCG) random variables all with unit variance, i.e., $\mathrm{E}\{\mathbf{d}_k \mathbf{d}_k^H\} = \mathbf{I}_{S_k}$, $\mathrm{E}\{\mathbf{d}_k \mathbf{d}_i^H\} = \mathbf{0}$, $\forall i \neq k$, and $\mathrm{E}\{\mathbf{d}_k \mathbf{n}_i^H\} = \mathbf{0}$. And $\mathbf{n}_k$ is a ZMCSCG random variable with covariance matrix $\mathbf{R}_{nk} \in \mathbb{C}^{M_k \times M_k}$.

To establish our new weighted sum MSE duality (for solving $\mathcal{P}1$), we model the virtual uplink channel system model as shown in Fig. 1.(b). This channel is modeled by introducing precoders $\{\mathbf{V}_k = [\mathbf{v}_{k1}, \cdots, \mathbf{v}_{kS_k}]\}_{k=1}^K$ and decoders $\{\mathbf{T}_k = [\mathbf{t}_{k1}, \cdots, \mathbf{t}_{kS_k}]\}_{k=1}^K$, where $\mathbf{v}_{ki} \in \mathbb{C}^{M_k \times 1}$ and $\mathbf{t}_{ki} \in \mathbb{C}^{N \times 1}, \forall k, i$. It is assumed that the $k$th user $i$th symbol ($d_{ki}$) is an i.i.d ZMCSCG random variable with variance $\zeta_{ki}$, i.e., $\mathrm{E}\{d_{ki} d_{ki}^H\} = \zeta_{ki}$, $\mathrm{E}\{d_{ki} d_{jm}^H\} = 0, \forall (j,m) \neq (k,i)$, and $\mathrm{E}\{\mathbf{d}_k \mathbf{n}_i^H\} = \mathbf{0}$. Moreover, $\mathbf{n}$ is a ZMCSCG random variable with diagonal covariance matrix $\mathbf{\Psi} = \mathrm{diag}(\psi_1, \cdots, \psi_N)$.

For the downlink system model of Fig. 1.a, the MSE between $d_{ki}$ and $\hat{d}_{ki}$ can be expressed as
$$\xi_{ki} = \mathrm{E}_\mathbf{d}\{(\hat{d}_{ki} - d_{ki})(\hat{d}_{ki} - d_{ki})^H\} \quad (2)$$
$$= \mathbf{w}_{ki}^H(\mathbf{H}_k^H \mathbf{B}\mathbf{B}^H \mathbf{H}_k + \mathbf{R}_{nk})\mathbf{w}_{ki} - 2\Re\{\mathbf{w}_{ki}^H \mathbf{H}_k^H \mathbf{b}_{ki}\} + 1.$$

When perfect CSI is available at the BS and MSs, the minimum MSE (MMSE) receiver of the $k$th user $i$th symbol is given as $\mathbf{w}_{ki} = (\mathbf{H}_k^H \mathbf{B}\mathbf{B}^H \mathbf{H}_k + \mathbf{R}_{nk})^{-1}\mathbf{H}_k^H \mathbf{b}_{ki}$. Plugging $\mathbf{w}_{ki}$ into (2), we get the MMSE of the $k$th user $i$th symbol as
$$\tilde{\tilde{\xi}}_{ki} = 1 - \mathbf{b}_{ki}^H \mathbf{H}_k (\mathbf{H}_k^H \mathbf{B}\mathbf{B}^H \mathbf{H}_k + \mathbf{R}_{nk})^{-1}\mathbf{H}_k^H \mathbf{b}_{ki}. \quad (3)$$

## III. Problem Formulation

Mathematically, the weighted sum rate maximization problem can be formulated as
$$\mathcal{P}1 : \max_{\{\mathbf{B}\}} \sum_{k=1}^K \sum_{i=1}^{S_k} \omega_{ki} R_{ki}, \text{ s.t } [\mathbf{B}\mathbf{B}^H]_{n,n} \leq \breve{p}_n, \forall n \quad (4)$$

where $R_{ki} = \log_2(\tilde{\tilde{\xi}}_{ki})^{-1}$ and $\omega_{ki} \geq 0$ are the achievable rate [3], [4] and weighting factor of the $k$th user $i$th symbol, respectively, and $\breve{p}_n$ is the maximum available power at the $n$th BS antenna. Without loss of generality, we assume that $\{0 < \omega_{ki} < 1, \forall i\}_{k=1}^K$. After several steps, the above problem can be equivalently formulated as [3] (see (16) of [3])
$$\min_{\{\tau_{ki},\nu_{ki},\mathbf{b}_{ki},\mathbf{w}_{ki},\forall i\}_{k=1}^K} \sum_{k=1}^K \sum_{i=1}^{S_k} \theta_{ki} \frac{1}{\tau_{ki}} \nu_{ki}^{\gamma_{ki}} + \sum_{k=1}^K \sum_{i=1}^{S_k} \eta_{ki} \xi_{ki},$$
$$\text{s.t } [\mathbf{B}\mathbf{B}^H]_{n,n} \leq \breve{p}_n, \prod_{k=1}^K \prod_{i=1}^{S_k} \nu_{ki} = 1, \nu_{ki} > 0, \tau_{ki} > 0 \quad (5)$$

where $\eta_{ki} = \tau_{ki}^{\mu_{ki}}$, $\gamma_{ki} = \frac{1}{1-\omega_{ki}}$, $\mu_{ki} = \frac{1}{\omega_{ki}} - 1$ and $\theta_{ki} = \omega_{ki} \mu_{ki}^{(1-\omega_{ki})}$. For fixed $\{\tau_{ki}, \nu_{ki}, \forall i\}_{k=1}^K$, the above optimization problem has the same mathematical structure as that of the downlink weighted sum MSE minimization problem. Therefore, by keeping $\{\tau_{ki}, \nu_{ki}, \forall i\}_{k=1}^K$ constant, $\{\mathbf{b}_{ki}, \mathbf{w}_{ki}, \forall i\}_{k=1}^K$ can be optimized by applying MSE downlink-uplink duality technique. For better exposition, let us summarize our duality based approach for solving $\mathcal{P}1$ as shown in **Algorithm I**.

In the following, we first establish our new weighted sum MSE downlink-uplink duality (i.e., to perform steps (2) and (4)). The current duality is able to transfer the weighted sum MSE from downlink to uplink and vice versa by using only one scaling factor[1]. Then, we employ the technique of [6] to ensure the power constraint of each BS antenna at step (4).

**Algorithm I**
Initialization: Set $\{\mathbf{B}_k\}_{k=1}^K$ such that $[\mathbf{B}\mathbf{B}]_{n,n} = p_n, \forall n$.
Update $\{\mathbf{W}_k\}_{k=1}^K$ using MMSE receiver approach, i.e.,
$$\mathbf{W}_k = (\mathbf{H}_k^H \mathbf{B}\mathbf{B}^H \mathbf{H}_k + \mathbf{R}_{nk})^{-1}\mathbf{H}_k^H \mathbf{B}_k, \forall k. \quad (6)$$
Update $\{\tau_{ki}, \nu_{ki}, \forall i\}_{k=1}^K$ by solving
$$\min_{\{\tau_{ki},\nu_{ki},\forall i\}_{k=1}^K} \sum_{k=1}^K \sum_{i=1}^{S_k} \theta_{ki} \frac{1}{\tau_{ki}} \nu_{ki}^{\gamma_{ki}} + \sum_{k=1}^K \sum_{i=1}^{S_k} \eta_{ki} \tilde{\tilde{\xi}}_{ki},$$
$$\prod_{k=1}^K \prod_{i=1}^{S_k} \nu_{ki} = 1, \nu_{ki} > 0, \tau_{ki} > 0 \quad \forall k, i. \quad (7)$$

**Repeat**
1) By the current $\{\mathbf{B}_k, \mathbf{W}_k\}_{k=1}^K$, compute the downlink weighted sum MSE (i.e., $\xi_w^{DL} = \sum_{k=1}^K \sum_{i=1}^{S_k} \eta_{ki} \xi_{ki}^{DL}$).
$$\xi_w^{DL} = \mathrm{tr}\{\boldsymbol{\eta}[\mathbf{W}^H \boldsymbol{\Gamma}^{DL} \mathbf{W} - 2\Re\{\mathbf{W}^H \mathbf{H}^H \mathbf{B}\} + \mathbf{I}]\} \quad (8)$$
where $\mathbf{H} = [\mathbf{H}_1, \cdots, \mathbf{H}_K]$, $\mathbf{W} = \mathrm{blkdiag}(\mathbf{W}_k, \forall k)$, $\boldsymbol{\eta} = \mathrm{diag}(\eta_{ki}, \forall k, i)$, $\mathbf{R}_n = \mathrm{blkdiag}(\mathbf{R}_{nk}, \forall k)$ and $\boldsymbol{\Gamma}^{DL} = \mathbf{H}^H \mathbf{B}\mathbf{B}^H \mathbf{H} + \mathbf{R}_n$.
**Uplink channel**
2) Transfer the total weighted sum MSE from downlink to uplink channel.
3) Update $\{\mathbf{T}_k\}_{k=1}^K$ using MMSE receiver technique.
**Downlink channel**
4) Transfer the total weighted sum MSE from uplink to downlink channel.
5) Update $\{\mathbf{W}_k\}_{k=1}^K$ by (6) and $\{\tau_{ki}, \nu_{ki}, \forall i\}_{k=1}^K$ by (7).
**Until** convergence.

## IV. Weighted Sum MSE Downlink-Uplink Duality

To establish this duality, we first compute the total weighted sum MSE in the uplink channel as
$$\xi_w^{UL} = \mathrm{tr}\{\boldsymbol{\lambda}(\mathbf{T}^H \boldsymbol{\Sigma} \mathbf{T} - 2\Re\{\boldsymbol{\zeta} \mathbf{T}^H \mathbf{H}\mathbf{V}\} + \boldsymbol{\zeta})\} \quad (9)$$
where $\boldsymbol{\zeta} = \mathrm{diag}(\{\zeta_{ki}, \forall i\}_{k=1}^K)$, $\boldsymbol{\lambda} = \mathrm{diag}(\{\lambda_{ki}, \forall i\}_{k=1}^K)$, $\boldsymbol{\Sigma} = \mathbf{H}\mathbf{V}\boldsymbol{\zeta}\mathbf{V}^H \mathbf{H}^H + \boldsymbol{\Psi}$, $\mathbf{V} = \mathrm{blkdiag}(\mathbf{V}_1, \cdots, \mathbf{V}_K)$ and $\mathbf{T} = [\mathbf{T}_1, \cdots, \mathbf{T}_K]$ with $\lambda_{ki}$ as the MSE weight of the $k$th user $i$th symbol. When $\boldsymbol{\zeta} = \boldsymbol{\lambda} = \mathbf{I}$, (9) turns to (5) of [6].

### A. Weighted sum MSE transfer (From downlink to uplink)

To apply this weighted sum MSE transfer for $\mathcal{P}1$, we set the uplink precoder, decoder, input variance and weights as
$$\mathbf{V} = \mathbf{W}, \ \mathbf{T} = \mathbf{B}, \ \boldsymbol{\zeta} = \boldsymbol{\eta}, \ \boldsymbol{\lambda} = \mathbf{I}. \quad (10)$$
Substituting (10) into (9) and equating $\xi_w^{DL} = \xi_w^{UL}$ gives
$$\tilde{\tau} = \sum_{n=1}^N \psi_n \tilde{p}_n = \tilde{\mathbf{p}}^T \boldsymbol{\psi} \quad (11)$$

---
[1]Note that for a total BS power constrained MSE-based problems, the approaches of [5] and [8] compute $S$ scaling factors to transfer the weighted sum MSE from uplink to downlink channel and vice versa. This computation requires $O(S^3)$ operations (see [5] and [8] for more details). Thus, the current duality is more computational efficient than that of the existing duality.

where $\tilde{\tau} = \sum_{k=1}^{K} \text{tr}\{\boldsymbol{\eta}_k \mathbf{W}_k^H \mathbf{R}_{nk} \mathbf{W}_k\}$, $\boldsymbol{\psi} = [\psi_1, \cdots, \psi_N]^T$, $\tilde{\mathbf{p}} = [\tilde{p}_1, \cdots, \tilde{p}_N]^T$, $\tilde{p}_n = \text{tr}\{\tilde{\mathbf{b}}_n^H \tilde{\mathbf{b}}_n\}$ and $\tilde{\mathbf{b}}_n^H$ is the $n$th row of $\mathbf{B}$. The above equation shows that by choosing any $\{\psi_n\}_{n=1}^N$ that satisfy (11), one can transfer the precoder/decoder pairs of the downlink channel to the corresponding decoder/precoder pairs of the uplink channel ensuring $\xi_w^{DL} = \xi_w^{UL_2}$, where $\xi_w^{UL_2}$ is the uplink weighted sum MSE at step 2 of **Algorithm I**. However, here $\{\psi_n\}_{n=1}^N$ should be selected in a way that $\mathcal{P}1$ can be solved by **Algorithm I**. To this end, we choose $\boldsymbol{\psi}$ as

$$\tilde{\tau} \geq \tilde{\mathbf{p}}^T \boldsymbol{\psi}. \tag{12}$$

By doing so, the uplink channel weighted sum MSE is less than the downlink channel weighted sum MSE (i.e., $\xi_w^{DL} \geq \xi_w^{UL_2}$). As will be clear later, to solve (5) with **Algorithm I**, $\boldsymbol{\psi}$ should be selected as in (12). This shows that step 2 of **Algorithm I** can be carried out with (10).

To perform step 3 of **Algorithm I**, we update $\mathbf{T}$ of (10) by using the uplink MMSE receiver which is expressed as

$$\mathbf{T} = \boldsymbol{\Sigma}^{-1} \mathbf{H} \mathbf{V} \boldsymbol{\zeta} = (\mathbf{H}\mathbf{W}\boldsymbol{\eta}\mathbf{W}^H\mathbf{H}^H + \boldsymbol{\Psi})^{-1}\mathbf{H}\mathbf{W}\boldsymbol{\eta} \tag{13}$$

where the second equality is obtained from (10) (i.e., $\mathbf{V}=\mathbf{W}$ and $\boldsymbol{\zeta} = \boldsymbol{\eta}$). The above expression shows that by choosing $\{\psi_n > 0\}_{n=1}^N$, we ensure that $(\mathbf{H}\mathbf{W}\boldsymbol{\eta}\mathbf{W}^H\mathbf{H}^H + \boldsymbol{\Psi})^{-1}$ exists. Next, we transfer the total weighted sum MSE from uplink to downlink channel (i.e., we perform step 4 of **Algorithm I**) and show that the latter MSE transfer ensures the power constraint of each BS antenna.

### B. Weighted sum MSE transfer (From uplink to downlink)

For a given total weighted sum MSE in the uplink channel with $\boldsymbol{\zeta} = \boldsymbol{\eta}$ and $\boldsymbol{\lambda} = \mathbf{I}$, we can achieve the same weighted sum MSE in the downlink channel (with the MSE weighting matrix $\boldsymbol{\eta}$) by using a nonzero scaling factor ($\beta$) which satisfies

$$\widetilde{\mathbf{B}} = \beta \mathbf{T}, \quad \widetilde{\mathbf{W}} = \mathbf{V}/\beta. \tag{14}$$

In this precoder/decoder transformation, we use the notations $\widetilde{\mathbf{B}}$ and $\widetilde{\mathbf{W}}$ to differentiate with the precoder and decoder matrices used in Section IV-A. By substituting (14) into $\xi_w^{DL}$ (with $\widetilde{\mathbf{B}}=\mathbf{B}$, $\widetilde{\mathbf{W}}=\mathbf{W}$) and then equating the resulting sum MSE with that of the uplink channel (9), $\beta$ can be determined as

$$\beta^2 = \frac{\tilde{\tau}}{\sum_{n=1}^N \psi_n \tilde{\mathbf{t}}_n^H \tilde{\mathbf{t}}_n} \tag{15}$$

where $\tilde{\mathbf{t}}_n^H$ is the $n$th row of the MMSE matrix $\mathbf{T}$ (13).

The power of each BS antenna in the downlink channel is thus given by

$$\widetilde{\mathbf{b}}_n^H \widetilde{\mathbf{b}}_n = \text{tr}\{\beta^2 \tilde{\mathbf{t}}_n^H \tilde{\mathbf{t}}_n\} = \frac{\tilde{\tau} \tilde{\mathbf{t}}_n^H \tilde{\mathbf{t}}_n}{\sum_{i=1}^N \psi_i \tilde{\mathbf{t}}_i^H \tilde{\mathbf{t}}_i} \leq \breve{p}_n, \forall n \tag{16}$$

where $\widetilde{\mathbf{b}}_n^H$ is the $n$th row of $\widetilde{\mathbf{B}}$. By multiplying both sides of (16) by $\psi_n$, we can rewrite the above expression as $\psi_n \geq f_n, \forall n$, where $f_n = \frac{\tilde{\tau}}{\breve{p}_n} \frac{\psi_n \tilde{\mathbf{t}}_n^H \tilde{\mathbf{t}}_n}{\sum_{i=1}^N \psi_i \tilde{\mathbf{t}}_i^H \tilde{\mathbf{t}}_i}$. Now, suppose that there exist $\{\psi_n > 0\}_{n=1}^N$ that satisfy

$$\psi_n = f_n, \ \forall n. \tag{17}$$

From (17), we get $\psi_n \breve{p}_n = f_n \breve{p}_n, \forall n$. It follows that

$$\sum_{n=1}^N \psi_n \breve{p}_n = \sum_{n=1}^N f_n \breve{p}_n = \tilde{\tau}. \tag{18}$$

The above expression shows that the solution of (17) satisfies (18). Moreover, as $\{\breve{p}_n \geq \tilde{p}_n\}_{n=1}^N$, the latter solution also satisfies (12). Therefore, for problem $\mathcal{P}1$, by choosing $\{\psi_n\}_{n=1}^N$ such that (17) is satisfied, step 4 can be performed.

Next, we show that there exists at least a set of feasible $\{\psi_n > 0\}_{n=1}^N$ that satisfy (17). In this regard, we consider the following Theorem [9].

*Theorem 1:* Let $(\mathbf{X}, \|.\|_2)$ be a complete metric space. We say that $F: \mathbf{X} \to \mathbf{X}$ is an almost contraction, if there exist $\kappa \in [0, 1)$ and $\chi \geq 0$ such that

$$\|F(\mathbf{x}) - F(\mathbf{y})\|_2 \leq \kappa \|\mathbf{x} - \mathbf{y}\|_2 + \chi \|\mathbf{y} - F(\mathbf{x})\|_2,$$
$$\forall \mathbf{x}, \mathbf{y} \in \mathbf{X}. \tag{19}$$

If $F$ satisfies (19), then the following holds true.
1) $\exists \mathbf{x} \in \mathbf{X} : \mathbf{x} = F(\mathbf{x})$.
2) For any initial $\mathbf{x}_0 \in \mathbf{X}$, the iteration $\mathbf{x}_{n+1} = F(\mathbf{x}_n)$ for $n = 0, 1, 2, \cdots$ converges to some $\mathbf{x}^\star \in \mathbf{X}$.
3) The solution $\mathbf{x}^\star$ is not necessarily unique.

*Proof:* See *Theorem 1.1* of [9].

By defining $F$ as $F(\boldsymbol{\psi}) \triangleq [f_1, f_2, \cdots, f_N]$ with $\{\psi_n \in [\epsilon > 0, (\tilde{\tau} - \epsilon \sum_{i=1, i \neq n}^N \breve{p}_i)/\breve{p}_n]\}_{n=1}^N$ (we use $\epsilon = \min(10^{-6}, \{\tilde{\tau}/\breve{p}_n\}_{n=1}^N)$ for simulation), it can be easily seen that $\|F(\boldsymbol{\psi}_1) - F(\boldsymbol{\psi}_2)\|_2$, $\|\boldsymbol{\psi}_1 - \boldsymbol{\psi}_2\|_2$ and $\|\boldsymbol{\psi}_2 - F(\boldsymbol{\psi}_1)\|_2$ are bounded, $\forall \boldsymbol{\psi}_1, \boldsymbol{\psi}_2 \in \boldsymbol{\psi}$. As a result, there exist $\kappa$ and $\chi$ that satisfy (19). Thus, $F(\boldsymbol{\psi})$ is an almost contraction, i.e.,

$$\boldsymbol{\psi}_{n+1} = F(\boldsymbol{\psi}_n), \text{ with } \boldsymbol{\psi}_0 = [\psi_{01}, \psi_{02}, \cdots, \psi_{0N}] \geq \epsilon \mathbf{1}_N,$$
$$\text{for } n = 0, 1, 2, \cdots \text{ converges} \tag{20}$$

where $\mathbf{1}_N$ is an $N$ length vector with each element equal to unity. As a result, there exists $\{\psi_n \geq \epsilon\}_{n=1}^N$ that satisfy (17) and its solution can be obtained using the above fixed point iterations. Once the appropriate $\{\psi_n\}_{n=1}^N$ is obtained, step 5 of **Algorithm I** is immediate and hence $\mathcal{P}1$ can be solved iteratively using this algorithm.

To increase the convergence speed of **Algorithm I**, $\{\nu_{ki}, \tau_{ki}, \forall i\}_{k=1}^K$ can be optimized jointly with the downlink powers of all symbols by GP approach like in [6]. Towards this end, we decompose the precoders and decoders of the downlink channel as

$$\mathbf{B}_k = \mathbf{G}_k \mathbf{P}_k^{1/2}, \quad \mathbf{W}_k = \mathbf{U}_k \boldsymbol{\alpha}_k \mathbf{P}_k^{-1/2}, \quad \forall k \tag{21}$$

where $\mathbf{P}_k = \text{diag}(p_{k1}, \cdots, p_{kS_k}) \in \Re^{S_k \times S_k}$, $\mathbf{G}_k = [\mathbf{g}_{k1} \cdots \mathbf{g}_{kS_k}] \in \mathbb{C}^{N \times S_k}$, $\mathbf{U}_k = [\mathbf{u}_{k1} \cdots \mathbf{u}_{kS_k}] \in \mathbb{C}^{M_k \times S_k}$ and $\boldsymbol{\alpha}_k = \text{diag}(\alpha_{k1}, \cdots, \alpha_{kS_k}) \in \Re^{S_k \times S_k}$ are the transmit power, unity norm transmit filter, unity norm receive filter and receiver scaling factor matrices of the $k$th user, respectively, i.e., $\{\mathbf{g}_{ki}^H \mathbf{g}_{ki} = \mathbf{u}_{ki}^H \mathbf{u}_{ki} = 1, \forall i\}_{k=1}^K$.

To simplify the GP formulation, we collect $\{\mathbf{P}_k, \mathbf{G}_k, \mathbf{U}_k$ and $\boldsymbol{\alpha}_k\}_{k=1}^K$ as $\mathbf{P} = \text{blkdiag}(\mathbf{P}_1, \cdots, \mathbf{P}_K) = \text{diag}(p_1, \cdots, p_S)$, $\mathbf{G} = [\mathbf{G}_1, \cdots, \mathbf{G}_K] = [\mathbf{g}_1, \cdots, \mathbf{g}_S]$, $\mathbf{U} = \text{blkdiag}(\mathbf{U}_1, \cdots, \mathbf{U}_K) = [\mathbf{u}_1, \cdots, \mathbf{u}_S]$ and $\boldsymbol{\alpha} = \text{blkdiag}(\boldsymbol{\alpha}_1, \cdots, \boldsymbol{\alpha}_K) = \text{diag}(\alpha_1, \cdots, \alpha_S)$. By defining $\xi =$

$[\xi_{1,1}^{DL}, \cdots, \xi_{K,S_K}^{DL}]^T = [\{\xi_l^{DL}\}_{l=1}^S]^T$, the $l$th symbol MSE in the downlink channel can be expressed as (refer [4] and [6] for the details about (21) and the above descriptions)

$$\xi_l^{DL} = p_l^{-1}[(\mathbf{D} + \boldsymbol{\alpha}^2 \boldsymbol{\Phi}^T)\mathbf{p}]_l + p_l^{-1}\alpha_l^2 \mathbf{u}_l^H \mathbf{R}_n \mathbf{u}_l \quad (22)$$

where

$$[\boldsymbol{\Phi}]_{l,j} = \begin{cases} |\mathbf{g}_l^H \mathbf{H} \mathbf{u}_j|^2, & \text{for } l \neq j \\ 0, & \text{for } l = j \end{cases} \quad (23)$$

$$[\mathbf{D}]_{l,l} = \alpha_l^2 |\mathbf{g}_l^H \mathbf{H} \mathbf{u}_l|^2 - 2\alpha_l \Re(\mathbf{u}_l^H \mathbf{H}^H \mathbf{g}_l) + 1, \quad (24)$$

$1 \le l(j) \le S$ and $\mathbf{p} = [p_1, \cdots, p_S]^T$. For convenience, let us also rewrite $\{\tau_{ki}, \nu_{ki}, \mu_{ki}, \omega_{ki}, \gamma_{ki}$ and $\theta_{ki}, \forall i\}_{k=1}^K$ as $\{\tau_l, \nu_l, \mu_l, \omega_l, \gamma_l$ and $\theta_l\}_{l=1}^S$, respectively. Now, for fixed $\mathbf{G}, \mathbf{U}$ and $\boldsymbol{\alpha}$, we can optimize $\{\tau_l, \nu_l, p_l\}_{l=1}^S$ by

$$\min_{\{\tau_l,\nu_l,p_l\}_{l=1}^S} \sum_{l=1}^S \theta_l \frac{1}{\tau_l} \nu_l^{\gamma_l} + \tau_l^{\mu_l} \xi_l^{DL}$$

$$\text{s.t } \boldsymbol{\varsigma}_n^T \mathbf{p} \le \breve{p}_n, \prod_{l=1}^S \nu_l = 1, \nu_l > 0, \tau_l > 0 \; \forall l. \quad (25)$$

where $\gamma_l = \frac{1}{1-\omega_l}$, $\mu_l = \frac{1}{\omega_l} - 1$, $\theta_l = \omega_l \mu_l^{(1-\omega_l)}$ and $\boldsymbol{\varsigma}_n^T \in \Re^{1 \times S} = \{|[\mathbf{G}]_{(n,i)}|^2\}_{i=1}^S$. The above problem is a GP for which global optimality is guaranteed [10]. The modified duality based algorithm for $\mathcal{P}1$ is thus summarized as follows.

**Algorithm II**

Initialization: Like in **Algorithm I** with $\mathbf{B}$ as the first $S$ columns of $\mathbf{H}$.

**Repeat**    Uplink channel
1) Set $\mathbf{V} = \mathbf{W}, \mathbf{T} = \mathbf{B}$ and compute $\{\psi_n\}_{n=1}^N$ by (17) with $\boldsymbol{\psi}_0 = \{\psi_{0n} = \tau / \sum_{i=1}^N \breve{p}_i\}_{n=1}^N$.
2) Update $\{\mathbf{T}_k\}_{k=1}^K$ with the MMSE approach (13).

     **Downlink channel**
3) Transfer the total weighted sum MSE from uplink to downlink channel using (14).
4) Decompose $\{\mathbf{B}_k$ and $\mathbf{W}_k\}_{k=1}^K$ as in (21). Then, formulate and solve the GP problem (25).
5) Using the updated $\{\mathbf{P}_k\}_{k=1}^K$, compute $\{\mathbf{B}_k\}_{k=1}^K$ by (21). Then, update $\{\mathbf{W}_k\}_{k=1}^K$ by the MMSE approach (6).

**Until** convergence.

**Convergence**: It can be shown that at each step of this algorithm, the objective function of (5) is non-increasing [6]. Thus, **Algorithm II** is guaranteed to converge.

It can be easily seen that **Algorithm II** can be applied to solve weighted sum Rate (MSE) based constrained with groups of BS antenna (total BS) power problem.

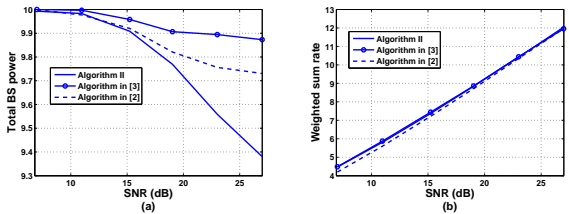

Fig. 2. Comparison of the proposed algorithm (**Algorithm II**) and the algorithms of [2] and [3] in terms of: (a) Total BS power. (b) Total weighted sum rate.

## V. SIMULATION RESULTS

For our simulation, we have used $K = 2$, $\{M_k = S_k = 2\}_{k=1}^K$ and $N = 4$. The channel between the BS and each MS consists of ZMCSCG entries with unit variance. The simulation results are averaged over 200 randomly chosen channel realizations. It is assumed that $\mathbf{R}_n = \sigma^2 \mathbf{I}$, $\{\breve{p}_n = 2.5\}_{n=1}^N$ and $\{\omega_l\}_{l=1}^4 = [0.4, 0.2, 0.6, 0.25]^T$. The Signal-to-Noise ratio (SNR) is defined as $P_{\text{sum}}/K\sigma^2$ and is controlled by varying $\sigma^2$, where $P_{\text{sum}}$ is the total BS power.

We compare the performance of our algorithm (**Algorithm II**) with that of [2] and [3]. The comparison is based on the total power utilized at the BS and total achieved weighted sum rate. As can be seen from Fig. 2.(a), the proposed algorithm utilizes less total BS power than that of [2] and [3]. Next, with the powers given in Fig. 2.(a), we plot the achieved weighted sum rate of these three algorithms which is shown in Fig. 2.(b). The latter figure shows that the proposed algorithm and the algorithm in [3] achieve the same weighted sum rate. From Fig 2.(a)-(b), one can observe that to achieve the same weighted sum rate, the proposed algorithm requires less total BS power than that of the algorithms in [2] and [3].

## VI. CONCLUSIONS

In this paper, we examine weighted sum rate maximization constrained with a per BS antenna power problem for downlink multiuser MIMO systems. To solve the problem, we propose novel downlink-uplink duality based iterative algorithm. Simulation results demonstrate that the proposed duality based iterative algorithm utilizes less total BS power compared to that of existing algorithms.